\documentclass{article}
\usepackage{amsmath}
\usepackage{amssymb}
\usepackage{amsfonts}

\newtheorem{theorem}{Theorem}

\newtheorem{theoreme}{Th\'{e}or\`{e}me}
\newtheorem{corllaire}{Corollaire}

\newtheorem{lemme}[theorem]{Lemme}

\newenvironment{proof}[1][Preuve]{\noindent\textbf{#1.} }{\ \rule{0.5em}{0.5em}}
\include {mak}
\begin{document}

\title{ Une identit\'{e} pour les polyn\^{o}mes de Bernoulli g\'{e}n\'{e}%
ralis\'{e}s}
\author{ Redha Chellal, Farid Bencherif \\
Faculty of Mathematics, USTHB, 16111 Algiers, Algeria.\\
rchellal@usthb.dz, fbencherif@usthb.dz, \and Mohamed Mehbali \\
Centre for Research Informed Teaching, LSBU, SE1 0AA London, UK. \\
mehbalim@lsbu.ac.uk\\
{\small \ } }
\maketitle

\begin{center}
\textbf{Résumé}
\end{center}
Dans cet article, nous \'{e}tablissons une identit\'{e} pour les polyn\^{o}%
mes g\'{e}n\'{e}ralis\'{e}s de Bernoulli. Nous en d\'{e}duisons des g\'{e}n%
\'{e}ralisations pour de nombreuses relations comportant les nombres ou les
polyn\^{o}mes de Bernoulli classiques. En particulier, nous g\'{e}n\'{e}%
ralisons une r\'{e}cente identit\'{e} de Gessel.

\begin{equation*}
\end{equation*}

\bigskip .


\section{Introduction}

Les polyn\^{o}mes de Bernoulli g\'{e}n\'{e}ralis\'{e}s $B_{n}^{\left( \alpha
\right) }\left( x\right) $ (\cite{rom}, p. 93) sont d\'{e}finis, pour $%
\alpha \in
\mathbb{C}
$, par :%
\begin{equation}
\sum_{n=0}^{\infty }B_{n}^{\left( \alpha \right) }\left( x\right) \frac{t^{n}%
}{n!}=\left( \frac{t}{e^{t}-1}\right) ^{\alpha }e^{tx}.  \label{cn1}
\end{equation}%
Les nombres de Bernoulli g\'{e}n\'{e}ralis\'{e}s $B_{n}^{\left( \alpha
\right) }$, les polyn\^{o}mes de Bernoulli classiques $B_{n}\left( x\right) $
et les nombres de Bernoulli $B_{n}$ sont respectivement d\'{e}finis par:%
\begin{equation}
B_{n}^{\left( \alpha \right) }=B_{n}^{\left( \alpha \right) }\left( 0\right)
\text{, \ }B_{n}\left( x\right) =B_{n}^{\left( 1\right) }\left( x\right)
\text{ et }B_{n}=B_{n}\left( 0\right) \text{.}  \label{t2}
\end{equation}%
Pour tous entiers naturels $n,\ell $ et $r$, consid\'{e}rons la somme $%
S_{n,l,r}^{\left( \alpha \right) }\left( x,y,z\right) $ d\'{e}finie par:%
\begin{eqnarray}
S_{n,l,r}^{\left( \alpha \right) }\left( x,y,z\right)
&=&\sum_{k=0}^{n+r}x^{n+r-k}\binom{n+r}{k}\binom{\ell +k+r}{r}B_{\ell
+k}^{\left( \alpha \right) }\left( y\right)  \notag \\
&&+(-1)^{\ell +n+r+1}\sum_{k=0}^{\ell +r}x^{\ell +r-k}\binom{\ell +r}{k}%
\binom{n+k+r}{r}B_{n+k}^{\left( \alpha \right) }\left( z\right) .  \label{t1}
\end{eqnarray}%
Il est bien connue \cite{ben1} que l'on a
\begin{equation}
S_{n,\ell \text{,},r}^{\left( 1\right) }\left( 1,0,0\right) =0.  \label{t3}
\end{equation}

En 2013, \`{a} l'aide du calcul ombral, Gessel \cite{bel} a g\'{e}n\'{e}ralis%
\'{e} ce r\'{e}sulat en prouvant la formule explicite suivante pour la somme
$S_{n,l,r}^{\left( 1\right) }\left( m,0,0\right) $ o\`{u} $m\geq 1$ est un
entier: \
\begin{equation}
S_{n,\ell \text{,},r}^{\left( 1\right) }\left( m,0,0\right) =\left(
r+1\right) \sum_{k=1}^{m-1}\sum_{j=0}^{r+1}\left( -1\right) ^{\ell +j-1}%
\binom{n+r}{j}\binom{\ell +r}{r+1-j}k^{\ell +j-1}(m-k)^{n+r-j}.  \label{t4}
\end{equation}

En consid\'{e}rant alors le cas particulier o\`{u} $n=\ell $ et o\`{u} $r$
est impair, Gessel \cite{bel} d\'{e}duit de la formule (\ref{t4}) la formule
explicite suivante:
\begin{equation}
\sum_{k=0}^{n+r}m^{n+r-k}\binom{n+r}{k}\binom{n+k+r}{r}B_{n+k}=\frac{1}{2}%
\left( r+1\right) \sum_{k=1}^{m-1}\sum_{j=0}^{r+1}\binom{n+r}{j}\binom{n+r}{%
r+1-j}k^{j+n-1}(k-m)^{n+r-j}.  \label{t5}
\end{equation}

Les identit\'{e}s (\ref{t3}), (\ref{t4}) et (\ref{t5}) sont en fait
l'aboutissement d'un long cheminement dont nous signalons certaines \'{e}%
tapes au troisi\`{e}me paragraphe. Notre principal r\'{e}sultat est une
expression simplifi\'{e}e pour la somme $S_{n,l,r}^{\left( \alpha \right)
}\left( x,y,z\right) $ dans le cas o\`{u} $x+y+z-\alpha $ est un entier
naturel. Nous l'\'{e}non\c{c}ons et le prouvons au quatri\`{e}me paragraphe.
Nous en d\'{e}duisons aussi de nombreux corollaires qui sont des g\'{e}n\'{e}%
ralisations d'identit\'{e}s connues. Au paragraphe suivant, nous rappelons
des propri\'{e}t\'{e}s des polyn\^{o}mes g\'{e}n\'{e}ralis\'{e}s de
Bernoulli et nous prouvons deux lemmes qui nous seront utiles pour prouver
notre r\'{e}sultat.

\section{Quelques propri\'{e}t\'{e}s des polyn\^{o}mes de Bernoulli g\'{e}n%
\'{e}ralis\'{e}s et lemmes}

Dans tout ce qui suit, nous convenons d'appeler op\'{e}rateur tout
endomorphisme du $%
\mathbb{C}
$-espace vectoriel $%
\mathbb{C}
\left[ x\right] $. Nous d\'{e}signons par $D$ l'op\'{e}rateur de d\'{e}%
rivation classique et par $I$ \ et $\Delta $ l'op\'{e}rateur identit\'{e} et
l'op\'{e}rateur de diff\'{e}rence finie d\'{e}finis respectivement par%
\begin{equation*}
I\left( x^{n}\right) =x^{n}\text{ et }\Delta \left( x^{n}\right) =\left(
x+1\right) ^{n}-x^{n}\text{, \ }n\in
\mathbb{N}
.
\end{equation*}

Les polyn\^{o}mes de Bernoulli g\'{e}n\'{e}ralis\'{e}s constituent une suite
de polyn\^{o}mes d'Appell. Ils v\'{e}rifient les propri\'{e}t\'{e}s
suivantes bien connues \cite{rom} et faciles \`{a} prouver:

\begin{equation}
B_{0}^{\left( \alpha \right) }\left( x\right) =1\text{ \ et }D\left(
B_{n}^{\left( \alpha \right) }\left( x\right) \right) =nB_{n-1}^{\left(
\alpha \right) }\left( x\right) \text{, }n\geq 1\text{,}  \label{f14}
\end{equation}%
\begin{equation}
B_{n}^{\left( \alpha \right) }\left( x+y\right) =\sum_{k=0}^{n}\binom{n}{k}%
y^{n-k}B_{k}^{\left( \alpha \right) }\left( x\right) ,  \label{f15}
\end{equation}

\begin{equation}
\Delta \left( B_{n}^{\left( \alpha \right) }\left( x\right) \right) =D\left(
B_{n}^{\left( \alpha -1\right) }\left( x\right) \right) ,  \label{f13}
\end{equation}%
\begin{equation}
B_{n}^{\left( \alpha \right) }\left( \alpha -x\right) =\left( -1\right)
^{n}B_{n}^{\left( \alpha \right) }\left( x\right) .  \label{f30}
\end{equation}

Pour tout $\gamma \in
\mathbb{C}
$, consid\'{e}rons l'endomorphisme $\Omega _{\gamma }$ de $%
\mathbb{C}
\left[ x\right] $ d\'{e}fini par:
\begin{equation*}
\Omega _{\gamma }\left( x^{n}\right) =B_{n}^{\left( \gamma \right) }\left(
x\right) \text{, \ }n\in
\mathbb{N}
.
\end{equation*}

\begin{lemme}
\label{lem4}Pour tout entier naturel $n$ et pour tous nombres complexes $%
\alpha $ et $\gamma $, on a:
\begin{equation}
\Omega _{\alpha }\left( \left( x+\gamma \right) ^{n}\right) =B_{n}^{\left(
\alpha \right) }\left( x+\gamma \right) \text{.}  \label{f12}
\end{equation}
\end{lemme}

\begin{proof}
Il suffit d'exploiter la propri\'{e}t\'{e} (\ref{f15}) pour $y=\gamma $. On
a
\begin{equation*}
\Omega _{\alpha }\left( \left( x+\gamma \right) ^{n}\right) =\Omega _{\alpha
}\left( \sum_{k=0}^{n}\binom{n}{k}\gamma ^{n-k}x^{k}\right) =\sum_{k=0}^{n}%
\binom{n}{k}y^{n-k}B_{k}^{\left( \alpha \right) }\left( x\right)
=B_{n}^{\left( \alpha \right) }\left( x+\gamma \right) .
\end{equation*}
\end{proof}

Les op\'{e}rateurs $\Omega _{\alpha }$, $\Delta $ et $D$ v\'{e}rifient de
remarquables relations donn\'{e}es dans le lemme suivant:

\begin{lemme}
\label{lem1}Pour tout $\alpha \in
\mathbb{C}
$, on a%
\begin{equation}
D\circ \Omega _{\alpha }=\Omega _{\alpha }\circ D.  \label{t9}
\end{equation}
\begin{equation}
\Omega _{\alpha }\circ \Delta =\Omega _{\alpha -1}\circ D.  \label{t6}
\end{equation}
\end{lemme}

\begin{proof}
\begin{enumerate}
\item Il suffit de v\'{e}rifier que l'on a $\left( D\circ \Omega _{\alpha
}\right) \left( x^{n}\right) =\left( \Omega _{\alpha }\circ D\right) \left(
x^{n}\right) $ pour tout entier $n\geq 0$. La propri\'{e}t\'{e} est \'{e}%
vidente pour $n=0$ car $B_{0}^{\left( \alpha \right) }\left( x\right) =1$ et
donc $\left( D\circ \Omega _{\alpha }\right) \left( x^{0}\right) =0=\left(
\Omega _{\alpha }\circ D\right) \left( x^{0}\right) $. Pour $n\geq 1$, on a
\begin{equation*}
\left( D\circ \Omega _{\alpha }\right) \left( x^{n}\right) =D\left(
B_{n}^{\left( \alpha \right) }\left( x\right) \right) =\Omega _{\alpha
}\left( nx^{n-1}\right) =\left( \Omega _{\alpha }\circ D\right) \left(
x^{n}\right) .
\end{equation*}

\item On a, en exploitant le lemme \ref{lem4} pour $\gamma =1$, les propri%
\'{e}t\'{e}s (\ref{f13}) et (\ref{t9}), pour tout entier $n\geq 0:$
\begin{eqnarray*}
\left( \Omega _{\alpha }\circ \Delta \right) \left( x^{n}\right) &=&\Omega
_{\alpha }\left( \left( x+1\right) ^{n}-x^{n}\right) \\
&=&B_{n}^{\left( \alpha \right) }\left( x+1\right) -B_{n}^{\left( \alpha
\right) }\left( x\right) \\
&=&\Delta \left( B_{n}^{\left( \alpha \right) }\left( x\right) \right) \\
&=&D\left( B_{n}^{\left( \alpha -1\right) }\left( x\right) \right) \\
&=&\left( D\circ \Omega _{\alpha -1}\right) \left( x^{n}\right) \\
&=&\left( \Omega _{\alpha -1}\circ D\right) \left( x^{n}\right) .
\end{eqnarray*}
\end{enumerate}
\end{proof}

Les propri\'{e}t\'{e}s suivantes des polyn\^{o}mes classiques de Bernoulli
sont bien connues \cite{abr}

\begin{equation}
\sum_{k=0}^{n}\binom{n}{k}B_{k}=B_{n}+\delta _{n,1},  \label{f20}
\end{equation}

\begin{equation}
\left( -1\right) ^{n}B_{n}=B_{n}+\delta _{n,1},  \label{f21}
\end{equation}%
\begin{equation}
B_{2n+1}=0\text{, }n\geq 1,  \label{f22}
\end{equation}%
o\`{u} $\delta _{i,j}$ est le symbole de Kronecker valant $1$ si $i=j$ et $0$
sinon.

Les propri\'{e}t\'{e}s suivantes sont bien connues et faciles \`{a} d\'{e}%
duire de la relation (\ref{cn1}) \'{e}crite pour $\alpha =1:$%
\begin{equation*}
\sum_{k=0}^{n}\binom{n}{k}B_{k}\left( x\right) =B_{n}\left( x\right) +x^{n},
\end{equation*}%
Ces relations permettent de calculer les premiers polyn\^{o}mes et nombres
de Bernoulli. Les polyn\^{o}mes et nombres de Bernoulli jouent un r\^{o}le
fondamental dans diff\'{e}rentes branches des math\'{e}matiques telles que
la combinatoire, la th\'{e}orie des nombres, l'analyse et la topologie. On
peut trouver dans \cite{dil} de tr\`{e}s nombreux articles concernant les
propri\'{e}t\'{e}s de ces nombres et polyn\^{o}mes.

\section{R\'{e}trospective}

Commen\c{c}ons tout d'abord par remarquer que les relations (\ref{t4}) et (%
\ref{t5}) prouv\'{e}es par Gessel peuvent se reformuler sous des formes \'{e}%
quivalentes qui vont nous \^{e}tre utiles pour cet expos\'{e}. En effet, en
changeant $k$ en $m-k$ au second membre de (\ref{t4}), on obtient l'identit%
\'{e} suivante \'{e}quivalente \`{a} celle donn\'{e}e par Gessel

\begin{equation}
S_{n,l,r}^{\left( 1\right) }\left( m,0,0\right) =\left( r+1\right)
\sum_{k=1}^{m-1}\sum_{j=0}^{r+1}\binom{n+r}{j}\binom{\ell +r}{r+1-j}%
k^{n+r-j}\left( k-m\right) ^{\ell +j-1}.  \label{tg4}
\end{equation}%
En remarquant que pour tous entiers naturels $r$ et $s$ tel que $r+s$ est
pair, on a la propri\'{e}t\'{e} suivante:
\begin{equation}
\sum_{k=1}^{m-1}\left( k^{r}(k-m)^{s}-k^{s}(k-m)^{r}\right) =0,  \label{rem1}
\end{equation}%
on peut aussi constater que la formule (\ref{t5}) \'{e}tablie par Gessel est
\'{e}quivalente \`{a} la formule donn\'{e}e dans \cite{aid1}:
\begin{equation}
S_{n,n,r}^{\left( 1\right) }\left( m,0,0\right) =\sum_{k=1}^{m-1}q_{k}(m,r,n)
\label{ges1}
\end{equation}%
o\`{u}
\begin{equation}
q_{k}(m,r,n)=\frac{r+1}{2}\binom{n+r}{\frac{r+1}{2}}^{2}\left( k\left(
m-k\right) \right) ^{n+\frac{r-1}{2}}+\left( r+1\right) \sum_{j=0}^{\frac{r-1%
}{2}}\binom{n+r}{j}\binom{n+r}{r+1-j}k^{j+n-1}(k-m)^{n+r-j}.  \label{tg5}
\end{equation}

Dans ce qui suit, nous allons constater que les identit\'{e}s de Gessel (\ref%
{t4}) et (\ref{t5}) ou leurs formulations \'{e}quivalentes (\ref{tg4}) et (%
\ref{ges1}) g\'{e}n\'{e}ralisent un grand nombre d'identit\'{e}s sur les
nombres de Bernoulli.

Ainsi la propri\'{e}t\'{e} bien connue suivante \cite{ben1} qu'on d\'{e}duit
ais\'{e}ment de (\ref{f20}) et (\ref{f21})
\begin{equation}
\sum_{k=0}^{n}\binom{n}{k}B_{k}=\left( -1\right) ^{n}B_{n},  \label{p1}
\end{equation}%
se traduit par la relation:%
\begin{equation*}
S_{n,0,0}^{\left( 1\right) }\left( 1,0,0\right) =0.
\end{equation*}%
La propri\'{e}t\'{e} (\ref{p1}) exprime que la suite $\left( \left(
-1\right) ^{n}B_{n}\right) _{n\geq 0}$ est autoduale (ou est une suite de C%
\'{e}saro ) \cite{ben1} ,\cite{zek2}. Il existe aussi de nombreuses formules
explicites pour les nombres de Bernoulli \cite{gou1}. Cependant, de nombreux
auteurs ont cherch\'{e} d'autres relations de r\'{e}currence \ permettant de
calculer plus facilement les nombres de Bernoulli. Ainsi, en 1880, Lucas
\cite{luc1} prouve \`{a} l'aide du calcul symbolique la relation suivante
qu'il qualifie d'importante:%
\begin{equation}
\sum_{k=0}^{n}\binom{n}{k}B_{\ell +k}+\left( -1\right) ^{\ell
+n+1}\sum_{k=0}^{\ell }\binom{\ell }{k}B_{n+k}=0,  \label{e1}
\end{equation}%
qui se traduit par%
\begin{equation}
S_{n,\ell ,0}^{\left( 1\right) }\left( 1,0,0\right) =0.  \label{e2}
\end{equation}%
et qui est donc un cas particulier de (\ref{t3}). Lucas remarque alors que dans
le cas particulier o\`{u} $\ell =n+1$, la formule (\ref{e2}) permet d'\'{e}%
crire:%
\begin{equation}
(n+1)S_{n,n+1,0}^{\left( 1\right) }\left( 1,0,0\right) =\sum_{k=0}^{n+1}%
\binom{n+1}{k}\left( n+k+1\right) B_{n+k}=0  \label{k5}
\end{equation}%
ce qui permet, compte tenu de (\ref{f22}), d'en d\'{e}duire que pour $n\geq 1$
\begin{equation}
\sum_{k=0}^{n}\binom{n+1}{k}\left( n+k+1\right) B_{n+k}=0
\label{k3}
\end{equation}%
Pour Lucas, la relation (\ref{k3}) qu'il rattache \`{a} des formules
obtenues en 1877 par Stern, pr\'{e}sente un grand int\'{e}r\^{e}t pour le
calul des nombres de Bernoulli. Elle permet d'exprimer le nombre de
Bernoulli $B_{2n}$ \`{a} l'aide des seuls nombres de Bernoulli $B_{j}$ pour $%
n\leq j\leq 2n-1$ alors qu'en appliquant la relation bien connue
\begin{equation}
B_{2n}=\frac{-1}{n+1}\sum_{k=0}^{2n-1}\binom{2n+1}{j}B_{j},
\end{equation}
ce m\^{e}me calcul de $B_{2n}$ n\'{e}c\'{e}ssiterait la connaissance des nombres de
Bernoulli $B_{j}$ pour $0\leq j\leq 2n-1$. En fait, la propri\'{e}t\'{e} (%
\ref{k3}) avait d\'{e}j\`{a} \'{e}t\'{e} aussi d\'{e}couverte en 1827 par
von Ettingshausen \cite{ett} avant d'\^{e}tre red\'{e}couverte en 1877 par L
Seidel \cite{sei}.

Lucas prouve la relation (\ref{e1}) qu'il qualifie d'importante, \`{a}
l'aide du calcul symbolique ou calcul ombral tel qu'on le retrouve d\'{e}%
velopp\'{e} dans les articles d'Agoh \cite{ago1} et de Gessel \cite{ges}.
Cette m\^{e}me relation a \'{e}t\'{e} ensuite l'objet d'une intense
recherche.

Apr\`{e}s Lucas, de nombreux auteurs demontrent (\ref{k3}) par diff\'{e}%
rentes m\'{e}thodes. Ainsi, Nielsen \cite{nie} la prouve aussi en 1923 et
Kaneko \  \cite{kan} en 1995. L'article de Kaneko contient deux preuves diff%
\'{e}rentes de cette relation, une preuve compliqu\'{e}e bas\'{e}e sur la th%
\'{e}orie des fractions continues appliqu\'{e}e aux s\'{e}ries formelles et
une seconde preuve beaucoup plus simple due \`{a} D. Zagier, bas\'{e}e sur
une transfomation involutive de suites. Elle a ensuite aussi \'{e}t\'{e}
prouv\'{e}e de nouveau en 2009 par Cigler (\cite{cig}) et en 2011 par Chang
et Ha \cite{cha} (Corollaire 1 (a)).

Ainsi, en 1971, Carlitz \cite{car} propose de d\'{e}duire la relation (\ref%
{e1}) uniquement de la relation (\ref{f20}). En 1972, Shannon \cite{sha}, en r\'{e}%
ponse au probl\`{e}me pos\'{e} par Carlitz, prouve cette relation par un
raisonnement par r\'{e}currence sur $m$ et $n$, en supposant toutefois prouv%
\'{e}e la relation (\ref{f22}).

La relation (\ref{e1}) est alors aussi prouv\'{e}e par Gessel en 2003 (lemme
7.2, \cite{ges}) par l'emploi du calcul ombral. En 2005, Vassilev et Missana
\cite{vas2} \ prouvent la relation suivante
\begin{equation*}
\sum_{k=0}^{n-1}\binom{n}{k}B_{\ell +k}+\left( -1\right) ^{\ell
+n+1}\sum_{k=0}^{\ell -1}\binom{\ell }{k}B_{n+k}=0,\text{ pour }n\geq 1\
\text{\ et }\ell \geq 1
\end{equation*}%
qui peut se traduire par%
\begin{equation*}
S_{n,\ell ,0}^{\left( 1\right) }\left( 1,0,0\right) -\left( 1-\left(
-1\right) ^{\ell +n}\right) B_{n+\ell }=0,\text{ pour }n\geq 1\  \text{\ et }%
\ell \geq 1
\end{equation*}%
et qui gr\^{a}ce \`{a} (\ref{f22}) et (\ref{p1}) est \'{e}quivalente \`{a} (%
\ref{e1}). Leur m\'{e}thode est bas\'{e}e sur l'\'{e}galit\'{e} des
coefficients de $x^{n}y^{m}$ et $x^{m}y^{n}$ dans un polyn\^{o}me sym\'{e}%
trique de deux variables $x$ et $y$ \ introduit dans \cite{vas1}.

La relation (\ref{e1}) est aussi prouv\'{e}e par par Chu et Magli \cite{chu}
en 2007 par l'exploitation de la relation (\ref{p1}), par Chen, L.H. Sun
\cite{che} en 2009 par une approche alg\'{e}brique informatique exploitant
une extension de l'algorithme de Zeilberger, par Bencherif et Garici \cite%
{ben1} en 2012 par l'exploitation de la relation (\ref{p1}), par Gould et
Quaintance en 2014 \cite{gou}, qui \'{e}tudie les propri\'{e}t\'{e}s de la
transformation binomiale d'une suite, par Prodinger en 2014 \cite{pro} qui
donne une preuve tr\`{e}s courte bas\'{e}e sur un emploi judicieux de
fonctions g\'{e}n\'{e}ratrices et par Neto en 2015 \cite{net} qui exploite
des propri\'{e}t\'{e}s de l'alg\`{e}bre de Zeon. N\'{e}to donne une g\'{e}n%
\'{e}ralisarion int\'{e}ressante mais comportant une erreur dans les
indices. L'\'{e}criture correcte \'{e}tant:
\begin{equation*}
\sum_{k=0}^{n}\alpha ^{n-k}\binom{n}{k}B_{\ell +k}^{\left( \alpha \right)
}+(-1)^{\ell +n+1}\sum_{k=0}^{\ell }\alpha ^{\ell -k}\binom{\ell }{k}%
B_{n+k}^{\left( \alpha \right) }=0.
\end{equation*}%
c'est \`{a} dire
\begin{equation*}
S_{n,l,0}^{\left( \alpha \right) }\left( \alpha ,0,0\right) =0\text{.}
\end{equation*}

En 2000, Agoh \cite{ago1} (relation (4.3)) prouve la g\'{e}n\'{e}ralisation
suivante de (\ref{e2}) par l'emploi de congruences et du calcul ombral%
\begin{equation}
S_{n,l,r}^{\left( 1\right) }\left( 1,0,0\right) =0.  \label{s3}
\end{equation}

En 2001, Momimaya \cite{mom} prouve par une m\'{e}thode p-adique une formule
\'{e}quivalente au cas particulier $r=1$ de (\ref{s3}).

En 2003, Gessel \cite{ges}, commence par g\'{e}n\'{e}raliser la relation de
Lucas-Carlitz (\ref{e1}) en prouvant que
\begin{equation}
S_{n,l,0}^{\left( 1\right) }\left( m,0,0\right) =\sum_{k=1}^{m-1}\left(
\left( n+\ell \right) k-mn\right) k^{n-1}\left( k-m\right) ^{\ell -1}.
\label{t230}
\end{equation}

Pour $m=1$, la relation (\ref{t230}) fournit effectivement la relation et
pour $\ell =n+1$, Gessel en d\'{e}duit la g\'{e}n\'{e}ralisation (\ref{t24})
de l'identit\'{e} de Kaneko.

En 2003, Gessel \cite{ges} prouve \`{a} l'aide du calcul ombral, (Theorem
7.3. de \cite{ges}) la g\'{e}n\'{e}ralisation suivante:

\begin{equation}
(n+1)S_{n,n+1,1}^{\left( 1\right) }\left( m,0,0\right)
:=\sum_{k=0}^{n+1}m^{n+1-k}\binom{n+1}{k}\left( n+k+1\right)
B_{n+k}=\sum_{k=1}^{m-1}p_{k}\left( m,1,n\right) .  \label{t24}
\end{equation}%
o\`{u}%
\begin{equation*}
p_{k}\left( m,1,n\right) =\left( n+1\right) ^{2}k^{n}\left( k-m\right)
^{n}+n\left( n+1\right) k^{n+1}\left( k-m\right) ^{n-1}\text{.}
\end{equation*}

En 2009, Chen et Sun \cite{che} s'int\'{e}resse au cas $r=3$. Ils prouvent
en exploitant une extension de l'algorithme de Zeilberger que:
\begin{equation}
S_{n,n,3}^{\left( 1\right) }\left( m,0,0\right) :=\sum_{k=0}^{n+3}m^{n+3-k}%
\binom{n+3}{k}\binom{n+k+3}{3}B_{n+k}=\sum_{k=1}^{m-1}p_{k}(m,3,n)
\label{c1}
\end{equation}%
\bigskip

o\`{u} l'expression $p_{k}(m,3,n)$ donn\'{e}e par Chen et Sun peut s'\'{e}%
crire comme suit:%
\begin{equation}
p_{k}(m,3,n)=q_{k}(m,3,n)+\binom{n+3}{3}\left( 3n+11\right) \left(
k^{n+2}\left( k-m\right) ^{n}-k^{n}\left( k-m\right) ^{n+2}\right)  \notag
\end{equation}%
En exploitant (\ref{rem1}) on constate que la formule (\ref{c1}) est \'{e}%
quivalente \`{a} la formule (\ref{ges1}) pour $r=3.$

\section{Principaux r\'{e}sultats}

Le th\'{e}or\`{e}me suivant fournit une expression simplifi\'{e}e pour la
somme $S_{n,\ell ,r}^{\left( \alpha \right) }\left( \lambda ,x\text{, }%
\alpha +s-\lambda -x\right) $ pour $\alpha ,\lambda \in
\mathbb{C}
$, $\ell $, $n$, $m$, $r$, $s\in
\mathbb{N}
$.

\begin{theoreme}
\label{theo1}Pour tous nombres complexes $\alpha $, $\lambda \in
\mathbb{C}
$ et pour tous entiers naturels $\ell $, $n$, $m$, $r$, $s$, on a:\bigskip
\begin{eqnarray}
&&\sum_{k=0}^{n+r}\lambda ^{n+r-k}\binom{n+r}{k}\binom{\ell +k+r}{r}B_{\ell
+k}^{\left( \alpha \right) }\left( x\right)  \notag \\
&&+(-1)^{\ell +n+r+1}\sum_{k=0}^{\ell +r}\lambda ^{\ell +r-k}\binom{\ell +r}{%
k}\binom{n+k+r}{r}B_{n+k}^{\left( \alpha \right) }\left( \alpha +s-\lambda
-x\right)  \notag \\
&=&\Omega _{\alpha -1}\left( \frac{D^{r+1}}{r!}\sum_{k=1}^{s}\left(
x-k\right) ^{\ell +r}(x+\lambda -k)^{n+r}\right) .  \label{le1}
\end{eqnarray}
\end{theoreme}

\begin{proof}
Consid\'{e}rons le polyn\^{o}me $P\left( x\right) $ d\'{e}fini par
\begin{equation}
P\left( x\right) =\sum_{k=1}^{s}\frac{D^{r}}{r!}\left( x-k\right) ^{\ell
+r}(x+\lambda -k)^{n+r}.  \label{t7}
\end{equation}%
On sait d'apr\`{e}s la relation \ref{t6} du lemme \ref{lem1} qu'on a
\begin{equation}
\Omega _{\alpha }\left( \Delta P\left( x\right) \right) =\Omega _{\alpha
-1}\left( DP\left( x\right) \right) .  \label{t70}
\end{equation}%
En posant
\begin{equation*}
P_{k}\left( x\right) =\frac{D^{r}}{r!}\left( \left( x-k\right) ^{\ell
+r}(x+\lambda -k)^{n+r}\right) \text{, \  \ }k\in
\mathbb{N}
,\quad\text{on a}
\end{equation*}%

\begin{eqnarray*}
\Delta P\left( x\right) &=&\sum_{k=1}^{s}\left( P_{k-1}\left( x\right)
-P_{k}\left( x\right) \right) \\
&=&P_{0}\left( x\right) -P_{s}\left( x\right) \\
&=&\frac{D^{r}}{r!}\left( x^{\ell +r}(x+\lambda )^{n+r}\right) -\frac{D^{r}}{%
r!}\left( \left( x-s\right) ^{\ell +r}(x+\lambda -s)^{n+r}\right) \\
&=&\frac{D^{r}}{r!}\left( \sum_{k=0}^{n+r}\lambda ^{n+r-k}\binom{n+r}{k}%
x^{\ell +k+r}\right) + \\
&&(-1)^{\ell +n+r+1}\sum_{k=0}^{\ell +r}\lambda ^{\ell +r-k}\binom{\ell +r}{k%
}\left( -1\right) ^{n+k}\left( x+\lambda -s\right) ^{n+k+r} \\
&=&\sum_{k=0}^{n+r}\lambda ^{n+r-k}\binom{n+r}{k}\binom{\ell +k+r}{r}x^{\ell
+k}+ \\
&&(-1)^{\ell +n+r+1}\sum_{k=0}^{\ell +r}\lambda ^{\ell +r-k}\binom{\ell +r}{k%
}\binom{n+k+r}{r}\left( -1\right) ^{n+k}\left( x+\lambda -s\right) ^{n+k}.
\end{eqnarray*}%
En calculant $\Omega _{\alpha }\left( \Delta P\left( x\right) \right) $ avec
cette derni\`{e}re expression de $\Delta P\left( x\right) $ et compte tenu
des propri\'{e}t\'{e}s (\ref{f12}) et (\ref{f30}), on constate que $\Omega
_{\alpha }\left( \Delta P\left( x\right) \right) $ est \'{e}gal au premier
membre de (\ref{le1}). Il imm\'{e}diat de constater que $\Omega _{\alpha
-1}\left( DP\left( x\right) \right) $ est \'{e}gal au second membre de $%
\Omega _{\alpha -1}\left( DP\left( x\right) \right) $. La relation (\ref{t70}%
) permet de conclure.
\end{proof}

Le th\'{e}or\`{e}me \ref{theo1} qui est notre principal r\'{e}sultat g\'{e}n%
\'{e}ralise un grand nombre d'identit\'{e}s comportant les nombres ou les
polyn\^{o}mes de Bernoulli classiques ou g\'{e}n\'{e}ralis\'{e}s. Dans les
applications de ce th\'{e}or\`{e}me, on sera souvent amen\'{e} \`{a}
exploiter la relation
\begin{equation*}
B_{n+k}^{\left( \alpha \right) }\left( \alpha +s-\lambda -x\right) =\left(
-1\right) ^{n+k}B_{n+k}^{\left( \alpha \right) }\left( x+\lambda -s\right) .
\end{equation*}

\subsection{\textbf{Application 1: g\'{e}n\'{e}ralisation de la formule de
Gessel aux polyn\^{o}mes classiques de Bernoulli.}}

Pour $\alpha =1$, on a $B_{n}^{\left( 1\right) }\left( x\right) =B_{n}\left(
x\right) $ et $\Omega _{\alpha -1}=\Omega _{0}=I$. En application le th\'{e}%
or\`{e}me \ref{theo1} dans ce cas particulier et en appliquant la formule
de Leibniz pour calculer le second membre de (\ref{le1}), on obtient le
corollaire suivant qui est une g\'{e}n\'{e}ralisation du th\'{e}or\`{e}me
que Gessel \cite{bel} a \'{e}tabli en 2013.

\label{coro1}\bigskip Pour tout nombre complexe $\lambda \in
\mathbb{C}
$ et pour tous entiers naturels $\ell $, $n$, $m$, $r$, $s$, on a:%
\begin{eqnarray*}
&&\sum_{k=0}^{n+r}\lambda ^{n+r-k}\binom{n+r}{k}\binom{\ell +k+r}{r}B_{\ell
+k}\left( x\right) \\
&&+(-1)^{\ell +n+r+1}\sum_{k=0}^{\ell +r}\lambda ^{\ell +r-k}\binom{\ell +r}{%
k}\binom{n+k+r}{r}B_{n+k}\left( 1+s-\lambda -x\right) \\
&=&\left( r+1\right) \sum_{k=1}^{s}\sum_{j=0}^{r+1}\binom{n+r}{j}\binom{\ell
+r}{r+1-j}\left( x-k\right) ^{\ell +j-1}(x+\lambda -k)^{n+r-j}.
\end{eqnarray*}%
\bigskip Pour un entier naturel $m$, avec $\lambda =m,\quad s=m-1$ et $x=0$ dans
la relation du corollaire \ref{coro1}, on retrouve effectivement la formule (%
\ref{t4}) due \`{a} Gessel. Signalons qu'en 2014, He \cite{he2} retrouve la
formule de Gessel par l'emploi de $q$-nombres et polyn\^{o}mes de Bernoulli.

\subsection{\textbf{Application 2: g\'{e}n\'{e}ralisation d'une identit\'{e}
de Nielsen et d'une identit\'{e} d'Agoh aux polyn\^{o}mes g\'{e}n\'{e}ralis%
\'{e}s de Bernoulli.}}

Pour \ $\lambda =m-2\beta $, $s=1$ et $x$ remplac\'{e} par $x+\beta $ dans
la relation (\ref{le1}), l'application du th\'{e}or\`{e}me \ref{theo1} nous
founit le corollaire suivant:

\begin{corllaire}
Pour tous nombres complexes $\alpha $, $\beta $ et pour tous entiers
naturels $\ell $, $n$, $m$ et $r$, on a:%
\begin{eqnarray}
&&\sum_{k=0}^{n+r}\left( m-2\beta \right) ^{n+r-k}\binom{n+r}{k}\binom{\ell
+k+r}{r}B_{\ell +k}^{\left( \alpha \right) }\left( x+\beta \right)  \notag \\
&&-\sum_{k=0}^{\ell +r}\left( -1\right) ^{r+\ell +k}\left( m-2\beta \right)
^{\ell +r-k}\binom{\ell +r}{k}\binom{n+k+r}{r}B_{n+k}^{\left( \alpha \right)
}\left( x-\beta +m-1\right)  \notag \\
&=&\Omega _{\alpha -1}\left( \frac{D^{r+1}}{r!}\left( x+\beta -1\right)
^{\ell +r}(x-\beta +m-1)^{n+r}\right) .    \label{f10}
\end{eqnarray}
\end{corllaire}

Pour $\alpha =1$ et $m=1$, le corollaire fournit l'identit\'{e} (rectifi\'{e}%
e et aux notations pr\`{e}s) due \`{a} Nielsen (relation (10) en page 182 de
\cite{nie}.

Pour $\beta =0$ et $m=1$, en constatant que le second membre de (\ref{f10})
peut s'\'{e}crire:
\begin{eqnarray*}
\frac{D^{r+1}}{r!}\left( \left( x-1\right) ^{\ell +r}x^{n+r}\right) &=&\frac{%
D^{r}}{r!}\left( \left( n+r\right) x^{n+r-1}\left( x-1\right) ^{\ell
+r}+\left( \ell +r\right) x^{n+r}\left( x-1\right) ^{\ell +r-1}\right) \\
&=&\left( n+r\right) \sum_{k=0}^{r}\binom{n+r-1}{k}\binom{\ell +r}{r-k}%
x^{n+r-k-1}\left( x-1\right) ^{\ell +k} \\
&&+\left( \ell +r\right) \sum_{k=0}^{r}\binom{\ell +r-1}{k}\binom{n+r}{r-k}%
x^{n+k}\left( x-1\right) ^{\ell +r-k-1},
\end{eqnarray*}%
on obtient, aux notations pr\`{e}s, la relation (3.4) (i) (rectifi\'{e}e)
d'Agoh \cite{ago1}.

\subsection{\textbf{Application 3: g\'{e}n\'{e}ralisation d'identit\'{e}s de
Neto, Sun, Chen, Wu, Sun et Pan}}

Pour $m=a$ et avec un changement de notations, le th\'{e}or\`{e}me \ref%
{theo1} nous fournit le corollaire suivant:

\begin{corllaire}
Si $x+y+z=\alpha $, alors on a
\begin{equation}
\left( -1\right) ^{n}\sum_{k=0}^{n+r}\binom{n+r}{k}\binom{\ell +k+r}{r}%
x^{n+r-k}B_{l+k}^{\left( \alpha \right) }\left( y\right) =\left( -1\right)
^{\ell +r}\sum_{k=0}^{\ell +r}\binom{\ell +r}{k}\binom{n+k+r}{r}x^{\ell
+r-k}B_{n+k}^{\left( \alpha \right) }\left( z\right)  \label{s1}
\end{equation}%
\begin{equation}
\sum_{k=0}^{n+r}\binom{n+r}{k}\binom{\ell +k+r}{r}x^{n+r-k}B_{l+k}^{\left(
\alpha \right) }\left( y\right) =\sum_{k=0}^{\ell +r}\binom{\ell +r}{k}%
\binom{n+k+r}{r}\left( -x\right) ^{\ell +r-k}B_{n+k}^{\left( \alpha \right)
}\left( x+y\right)  \label{s2}
\end{equation}%
Si $x+y+z=s+1$ o\`{u} $s$ est un entier, alors on a%
\begin{equation}
S_{n,l,r}^{\left( 1\right) }\left( x,y,z\right) =\left( r+1\right)
\sum_{k=1}^{s}\sum_{j=0}^{r+1}\binom{n+r}{j}\binom{\ell +r}{r+1-j}\left(
y-k\right) ^{\ell +j-1}(x+y-k)^{n+r-j}  \label{s4}
\end{equation}
\end{corllaire}

Les relations (\ref{s1}) et (\ref{s2}) sont \'{e}quivalentes. Cela r\'{e}%
sulte du fait que $$B_{n+k}^{\left( \alpha \right) }\left( z\right)
=B_{n+k}^{\left( \alpha \right) }\left( \alpha -\left( x+y\right) \right)
=\left( -1\right) ^{n+k}B_{n+k}\left( x+y\right) $$. La relation (\ref{s1})
exprime que l'on a $S_{n,l,r}^{\left( \alpha \right) }\left( x,y,z\right) =0$
\ quand $x+y+z=\alpha .$

Pour $x=\alpha $, $y=z=0$ et $r=1$ dans (\ref{s1}), on retrouve l'identit%
\'{e} (rectifi\'{e}e) de Neto \cite{net}.

Pour $r=0$ et $\alpha =1$ dans (\ref{s1}), on retrouve la relation (1.15)
obtenue en 2003 par Sun \cite{sun}. Cette derni\`{e}re relation a aussi \'{e}%
t\'{e} prouv\'{e}e par Chen et Sun en 2009 (th\'{e}or\`{e}me 5.1) \cite{che}.

Pour $\alpha =1$ dans (\ref{s2}), on retrouve la relation (5.7) obtenue par
Chen en 2007 \cite{chen}.

Pour $\alpha =1$, $x=1$, $z=-y$ et $r=0$, la relation (\ref{s1}) permet
d'obtenir la relation $\left( 6\right) $ du th\'{e}or\`{e}me 2 de Wu, Sun et
Pan \cite{wu}.

On d\'{e}duit ais\'{e}ment de ce dernier corollaire l'identit\'{e} donn\'{e}%
e dans le corollaire suivant:

\begin{corllaire}
Si $x+y+z=\alpha $, alors on a pour $r\geq 1$, on a%
\begin{eqnarray*}
&&\left( -1\right) ^{n}\sum_{k=0}^{n+r-1}\binom{n+r}{k}\binom{\ell +k+r}{r}%
x^{n+r-k}B_{l+k}^{\left( \alpha \right) }\left( y\right) +\left( -1\right)
^{\ell +r+1}\sum_{k=0}^{\ell +r-1}\binom{\ell +r}{k}\binom{n+k+r}{r}x^{\ell
+r-k}B_{n+k}^{\left( \alpha \right) }\left( z\right) \\
&=&\left( -1\right) ^{n}\binom{\ell +n+2r}{r}\left( B_{n+\ell +1}^{\left(
\alpha \right) }\left( x+y)\right) -B_{^{{}}n+\ell +1}^{\left( \alpha
\right) }\left( y\right) \right) .
\end{eqnarray*}
\end{corllaire}

Pour $r=1$, $\alpha =1$, on retrouve la relation (1.16) obtenue en 2003 par
Sun \cite{sun}.

Pour $\alpha =1$, $x=1,y=t$, $z=-t$, on d\'{e}duit de ce corollaire le
corollaire suivant:

\begin{corllaire}
Pour $r\geq 1$, on a%
\begin{eqnarray*}
&&\left( -1\right) ^{n}\sum_{k=0}^{n+r-1}\binom{n+r}{k}\binom{\ell +k+r}{r}%
B_{l+k}\left( t\right) +\left( -1\right) ^{\ell +r+1}\sum_{k=0}^{\ell +r-1}%
\binom{\ell +r}{k}\binom{n+k+r}{r}B_{n+k}\left( -t\right) \\
&=&\left( -1\right) ^{n}\binom{n+\ell +2r}{r}\left( n+\ell +1\right)
t^{n+\ell }.
\end{eqnarray*}
\end{corllaire}

Pour $r=1$, on obtient l'identit\'{e} $\left( 8\right) $ du th\'{e}or\`{e}me
2 de Wu, Sun et Pan \cite{wu}.

Consid\'{e}rons la relation (\ref{s1}) pour $r$ impair, $\ell =n$ et $y=z=t$%
, on a alors $x=\alpha -2t$ et cette relation devient%
\begin{equation*}
\sum_{k=0}^{n+r}\binom{n+r}{k}\binom{n+k+r}{r}\left( \alpha -2t\right)
^{n+r-k}B_{n+k}^{\left( \alpha \right) }\left( t\right) =0.
\end{equation*}%
On en d\'{e}duit le corollaire suivant:

\begin{corllaire}
Pour tout entier naturel $n$ et pour tout entier naturel $r$ impair, on a
\begin{equation}
\sum_{k=0}^{n+r-1}\binom{n+r}{k}\binom{n+k+r}{r}\left( \alpha -2t\right)
^{n+r-k}B_{n+k}^{\left( \alpha \right) }\left( t\right) =-\binom{2n+2r}{r}%
B_{2n+r}^{\left( \alpha \right) }\left( t\right) \text{.}  \label{s20}
\end{equation}
\end{corllaire}

Pour $\alpha =1$ et $r=1$, la relation (\ref{s20}) fournit la relation
(1.18) de Sun \cite{sun}.

\subsection{\textbf{Application 4: g\'{e}n\'{e}ralisation d'une identit\'{e}
de Chang et Ha}}

\begin{corllaire}
Pour tout entier naturel $n$ et pour tout entier naturel $r$ impair, on a
\begin{equation}
\sum_{k=0}^{n+r}\binom{n+r}{k}\binom{n+k+r}{r}\frac{B_{n+k}\left( x\right) }{%
2^{k}\left( 1-x\right) ^{n+k-1}}=\left( -1\right) ^{n+\frac{r+1}{2}}\frac{r+1%
}{2^{n+r}}\binom{n+r}{\frac{r+1}{2}}.  \label{cor1}
\end{equation}%
\begin{equation}
\sum_{k=n}^{2n+r}\binom{n+r}{k-n}\binom{k+r}{r}\frac{B_{k}}{2^{k}}=\frac{%
\left( -1\right) ^{n+\frac{r-1}{2}}\left( r+1\right) }{2^{2n+r+1}}\binom{n+r%
}{\frac{r+1}{2}}.  \label{fi2}
\end{equation}
\end{corllaire}

La relation (\ref{fi2}) est \'{e}quivalente \`{a} la relation (\ref{cor1}).
Elle s'obtient par un simple changement d'indices.

\begin{proof}
Pour $\alpha =1$, $m=2$, $a=1$ et $x=0$, la relation (\ref{le1}) du th\'{e}or%
\`{e}me \ref{theo1} s'\'{e}crit%
\begin{eqnarray*}
2\sum_{k=0}^{n+r}\left( 2-2\beta \right) ^{n+r-k}\binom{n+r}{k}\binom{n+k+r}{%
r}B_{n+k}\left( \beta \right) &=&\left( r+1\right) \left[ \frac{D^{r+1}}{%
\left( r+1\right) !}\left( x^{2}-\left( \beta -1\right) ^{2}\right) ^{n+r}%
\right] _{x=0} \\
&=&\left( r+1\right) \binom{n+r}{\frac{r+1}{2}}\left( -\left( \beta
-1\right) ^{2}\right) ^{n+\frac{r-1}{2}}.
\end{eqnarray*}%
En rempla\c{c}ant $\beta $ par $x,$ le r\'{e}sultat du corollaire en r\'{e}sulte.
\end{proof}

Dans le cas particulier o\`{u} $r=1$ dans (\ref{fi2}) et $n\geq 1$, on
obtient ainsi la relation $(b)$ du corollaire 1 de Chang et Ha \cite{cha}.

\end{document}